\documentclass[10pt]{amsart}
\usepackage{amssymb,amsthm,amsfonts,amstext}
\usepackage{amsmath}
\usepackage{graphicx}
\usepackage[american]{babel}
\usepackage[ansinew]{inputenc}
\makeatletter \@addtoreset{equation}{section} \makeatother

\renewcommand\thetable{\thesection.\@arabic\c@table}
\newtheorem{theorem}{Theorem}[section]
\newtheorem{lemma}[theorem]{Lemma}
\newtheorem{proposition}[theorem]{Proposition}
\newtheorem{corollary}[theorem]{Corollary}

\newtheorem{remark}{Remark}[section]

\usepackage{a4wide}

\title[Equilibrium Fluctuations for the TAZRP]{Equilibrium Fluctuations for the\\ Totally Asymmetric Zero Range process}
\date{\today}
\begin{document}
\author{Patrícia Gonçalves}
\begin{abstract}
We prove a Central Limit Theorem for the empirical measure in the one-dimensional Totally Asymmetric Zero-Range Process in the hyperbolic
scaling $N$, starting from the equilibrium measure $\nu_{\rho}$. We also show that when taking the direction of the characteristics, the limit
density fluctuation field does not evolve in time until $N^{4/3}$, which implies the current across the characteristics to vanish in this longer
time scale.
\end{abstract}
\renewcommand{\subjclassname}{\textup{2000} Mathematics Subject Classification}
\begin{thanks} {The author was supported by
    F.C.T. (Portugal) with the grant /SFRH/ BPD/ 39991/ 2007.}
\end{thanks}
\address{Centro de Matemática da Universidade do Minho, Campus de Gualtar, 4710-057 Braga, Portugal}
\email{patg@math.uminho.pt}

\maketitle
\section{Introduction}
In this paper we study the Totally Asymmetric Zero-Range process (TAZRP) in $\mathbb{Z}$. In this process, if particles are present at a site
$x$, after a mean one exponential time, one of them jumps to $x+1$ at rate $1$, independently of particles on other sites. This is a Markov
process $\eta_{\cdot}$ with space state $\mathbb{N}^{\mathbb{Z}}$ and configurations are denoted by $\eta$, so that for a site $x$, $\eta(x)$
represents the number of particles at that site. For each density of particles $\rho$ there exists an invariant measure denoted by $\nu_{\rho}$,
which is translation invariant and such that $E_{\nu_{\rho}}[\eta(0)]=\rho$.

Since the work of Rezakhanlou in \cite{R.}, it is known that for the TAZRP the macroscopic particle density profile in the Euler scaling of
time, evolves according to the hyperbolic conservation law $\partial_{t}\rho(t,u)+\nabla \phi(\rho(t,u))=0$, where
$\phi(\rho)=\frac{\rho}{1+\rho}$. Since $\phi$ is differentiable, last equation can also be written as
$\partial_{t}\rho(t,u)+\phi'(\rho(t,u))\nabla\rho(t,u)=0$. This result is a Law of Large Numbers for the empirical measure associated to this
process starting from a general set of initial measures associated to a profile $\rho_{0}$, see \cite{R.} for details. If one wants to go
further and show a Central Limit Theorem (C.L.T.) for the empirical measure starting from the equilibrium state $\nu_{\rho}$, one has to
consider the density fluctuation field as defined below, see (\ref{eq:densfieldinz}).

Taking the hyperbolic time scale, the limit density field at time $t$ is just a translation of the initial density field. The translation or
velocity of the system is given by $\phi'(\rho)=\frac{1}{(1+\rho)^2}$ which is the characteristics speed. If we consider the particle system
moving in a reference frame with this constant velocity, then the limit field does not evolve in time and one is forced to consider a longer
time scale. Following the same approach as in \cite{G.} we can accomplish the result up to the time scale $N^{4/3}$, ie the limit density field
does not evolve in time until this time scale. Using this approach, the main difficulty in proving the C.L.T. for the empirical measure is the
Boltzmann-Gibbs Principle, which we can handle by using a multi-scale argument as done for the ASEP in \cite{G.}, but in this case there are
some extra computations to overcome the large space state. This result implies that the flux of particles through the characteristics speed
vanishes in this longer time scale. In fact, it was recently proved by \cite{B.K.} that the variance of the current across a characteristic is
of order $t^{2/3}$ and this translates by saying that in fact our result should hold till the time scale $N^{3/2}$. These results should be
valid for more general systems than TAZRP or TASEP (see \cite{G.}), but for systems with one conserved quantity and hyperbolic conservation law.
This is a step for showing this universality behavior.

This paper is a natural continuation of \cite{G.} and the multi-scale argument seems to be robust enough to be able to generalize it to other
models and to achieve the conjectured sharp time scale $N^{3/2}$, this is subject to future work.

We remark that all the results presented here, also hold for a more general Zero-Range process, namely one could take a Zero-Range dynamics in
which the jump rate from $x$ to $x+1$ is given by $g(\eta(x))$, with $g$ satisfying conditions of definition 3.1 of Chap.2 of \cite{K.L.}. We
could also consider a partial asymmetric process, in which a particle jumps from $x$ to $x+1$ at rate $pg(\eta(x))$ and from $x$ to $x-1$ at
rate $qg(\eta(x))$, where $p+q=1$ and $p\neq{1/2}$ and with $g$ as general as above. The results are valid for these more general processes but
in order to keep the presentation simple we state and prove them for the TAZRP.

An outline of the article follows. In the second section we introduce the notation and state the main results. In the third section we show the
Central Limit Theorem for the empirical measure in the hyperbolic time scale and the Central Limit Theorem for the current over a fixed bond. In
the fourth section we use the same approach as in \cite{G.} to prove the Central Limit Theorem for the empirical measure on a longer time scale
and the proof of the Boltzmann-Gibbs Principle is postponed to the fifth section.

\section{Statement of Results}

The generator of the one-dimensional TAZRP is given on local functions $f:\mathbb{N}^{\mathbb{Z}}\rightarrow{\mathbb{R}}$ by
\begin{equation*}
\mathcal{L}f(\eta)=\sum_{x\in{\mathbb{Z}}}1_{\{\eta(x)\geq{1}\}}[f(\eta^{x,x+1})-f(\eta)], \label{generatorTASEP}
\end{equation*}
where
\[\eta^{x,x+1}(z)=
\begin{cases}
\eta(z), & \mbox{if $z\neq{x,x+1}$}\\
\eta(x)-1, & \mbox{if $z=x$}\\
\eta(x+1)+1, & \mbox{if $z=x+1$}
\end{cases}.
\]
\\
In order to keep notation the more general as we can, we denote by $g(\eta(x))$ the function $1_{\{\eta(x)\geq{1}\}}$, which denotes the jump
rate of a particle to leave the site $x$.

The description of the process is the following. At each site, one can have any integer number of particles and for a site $x$ after an
exponential time of rate one, one of the particles at that site, jumps to the neighboring right site $x+1$, at rate 1. Initially, place the
particles according to a Geometric product measure in $\mathbb{N}^{\mathbb{Z}}$ of parameter $\rho$, denoted by $\nu_{\rho}$, which is an
invariant measure for the process.

Since the work of Rezakhanlou \cite{R.} it is known that taking the TAZRP in the Euler time scaling and starting from general initial measures
associated to an initial profile $\rho_{0}$ (for details we refer the reader to \cite{R.}), one gets in the hydrodynamic limit to the hyperbolic
conservation law:
\begin{equation*}
\partial_{t}\rho(t,u)+\nabla \phi(\rho(t,u))=0,
\end{equation*}
where the flux is given by $\phi(\rho)=\frac{\rho}{1+\rho}$.

Fixed a configuration $\eta$, let $\pi^{N}(\eta,du)$ denote the empirical measure given by
\begin{equation*}
\pi^{N}(\eta,du)=\frac{1}{N}\sum_{x\in\mathbb{Z}}\eta(x)\delta_{\frac{x}{N}}(du)
\end{equation*}
where $\delta_{u}$ denotes the Dirac measure at $u$ and let $\pi_{t}^{N}(\eta,du)=\pi^{N}(\eta_{t},du)$.

In order to state the C.L.T. for the empirical measure we need to define a suitable set of test functions. For an integer $k\geq{0}$, denote by
$\mathcal{H}_{k}$ the Hilbert space induced by $\mathcal{S}(\mathbb{R})$ (the Schwartz space) and the scalar product $<f,g>_{k}=<f,K_{0}^{k}g>$,
where $<\cdot,\cdot>$ denotes the inner product of $L^{2}(\mathbb{R})$ and $K_{0}$ is the operator $K_{0}=x^{2}-\Delta$. Denote by
$\mathcal{H}_{-k}$ the dual of $\mathcal{H}_{k}$, relatively to the inner product of $L^{2}(\mathbb{R})$ .

Fix $\rho$ and an integer $k$. Denote by $\mathcal{Y}_{.}^{N}$ the linear functional acting on functions $H\in \mathcal{S}(\mathbb{R})$ as
\begin{equation*}
\mathcal{Y}_{t}^{N}(H)=\sqrt{N}\Big[<H,\pi^{N}_{t}(\eta,du)>-\mathbb{E}_{\nu_{\rho}}<H,\pi^{N}_{t}(\eta,du)>\Big]
\end{equation*}
\begin{equation}
=\frac{1}{\sqrt{N}}\sum_{x\in{\mathbb{Z}}}H\Big(\frac{x}{N}\Big)(\eta_{tN}(x)-\rho). \label{eq:densfieldinz}
\end{equation}
where $<H,\pi^{N}_{t}(\eta,du)>$ denotes the integral of a test function $H$ wrt to the measure $\pi^{N}_{t}(\eta,du)$. Throughout the article
the functional above is mentioned as the density fluctuation field of the process. Denote by $D(\mathbb{R}^{+},\mathcal{H}_{-k})$ (resp.
$C(\mathbb{R}^{+},\mathcal{H}_{-k})$) the space of $H_{-k}$-valued functions, right continuous with left limits (resp. continuous), endowed with
the uniform weak topology, by $Q_{N}$ the probability measure on $D(\mathbb{R}^{+},\mathcal{H}_{-k})$ induced by $\mathcal{Y}^{N}_{.}$ and
$\nu_{\rho}$. Consider $\mathbb{P}^{N}_{\nu_{\rho}}=\mathbb{P}_{\nu_{\rho}}$ the p.m. on $D(\mathbb{R}^{+},\{0,1\}^{\mathbb{Z}})$ induced by
$\nu_{\rho}$ and $\eta_{t}$ speeded up by $N$ and denote by $\mathbb{E}_{\nu_{\rho}}$ the expectation with respect to $\mathbb{P}_{\nu_{\rho}}$.

\begin{theorem} \label{th:flu1}
Fix an integer $k>2$. Denote by $Q$ be the probability measure on $C(\mathbb{R}^{+},\mathcal{H}_{-k})$ corresponding to a stationary Gaussian
process with mean $0$ and covariance given by
\begin{equation*}
E_{Q}[\mathcal{Y}_{t}(H)\mathcal{Y}_{s}(G)]=\chi(\rho)\int_{\mathbb{R}}H(u+\phi'(\rho)(t-s))G(u)du \label{eq:covar}
\end{equation*}
for every $0\leq{s}\leq{t}$ and $H$, $G$ in $\mathcal{H}_{k}$. Here $\chi(\rho)= \textbf{Var}(\eta(0),\nu_{\rho})$. Then, $(Q_{N})_{N\geq{1}}$
converges weakly to $Q$.
\end{theorem}

Last result holds for the TAZRP evolving in any $\mathbb{Z}^{d}$ and for the other more general processes as pointed out in the introduction.

An easy consequence of last result is the derivation of the C.L.T. for the current over a fixed bond, see \cite{J.L.}. For a site $x$, denote by
$J^N_{x,x+1}(tN)$ the current through the bond $[x,x+1]$, defined as the total number of jumps from $x$ to $x+1$ during the time interval
$[0,tN]$.

\begin{theorem} \label{th:cltcurrent}
Fix $x\in{\mathbb{Z}}$, $t\geq{0}$ and let
\begin{equation*}
Z_{t}^{N}=\frac{1}{\sqrt{N}}\Big\{J^N_{x,x+1}(tN)-\mathbb{E}_{\nu_{\rho}}[J^N_{x,x+1}(tN)]\Big\}.
\end{equation*}
Then, under $\mathbb{P}_{\nu_{\rho}}$,
\begin{equation*}
\frac{Z_{t}^{N}}{\sqrt{\chi(\rho)\phi'(\rho)}}\xrightarrow[N\rightarrow{+\infty}]\,B_{t}
\end{equation*}
weakly, where $B_{t}$ denotes the standard Brownian motion.
\end{theorem}

In the hyperbolic scaling, the limit density fluctuation field at time $t$ is a translation of the initial one and removing from the system this
translation velocity, it does not evolve in time and one is forced to go beyond the hydrodynamic time scale. With this in mind, let
$\eta_{\cdot}$ be evolving in the time scale $N^{1+\gamma}$, with $\gamma>0$, fix $\rho$ and remove the translation velocity by redefining the
density fluctuation field on $H\in{\mathcal{S}(\mathbb{R})}$ by:
\begin{equation*}
\mathcal{Y}_{t}^{N,\gamma}(H)=\frac{1}{\sqrt{N}}\sum_{x\in{\mathbb{Z}}}H\Big(\frac{x-\phi'(\rho)tN^{1+\gamma}}{N}\Big)(\eta_{t{N}^{1+\gamma}}(x)-\rho).
\label{eq:densfieldlongscale}
\end{equation*}

As above, let $Q^\gamma_{N}$ be the probability measure on $D(\mathbb{R}^{+},\mathcal{H}_{-k})$ induced by $\mathcal{Y}^{N,\gamma}_{.}$ and
$\nu_{\rho}$, let $\mathbb{P}^{N,\gamma}_{\nu_{\rho}}=\mathbb{P}^{\gamma}_{\nu_{\rho}}$ be the p.m. on $D(\mathbb{R}^{+},\{0,1\}^{\mathbb{Z}})$
induced by $\nu_{\rho}$ and $\eta_{t}$ speeded up by $N^{1+\gamma}$ and denote by $\mathbb{E}_{\nu_{\rho}}^{\gamma}$ expectation with respect to
$\mathbb{P}_{\nu_{\rho}}^{\gamma}$. Then

\begin{theorem}
Fix an integer $k>1$ and $\gamma<1/3$. Let $Q$ be the probability measure on $C(\mathbb{R}^{+},\mathcal{H}_{-k})$ corresponding to a stationary
Gaussian process with mean $0$ and covariance given by
\begin{equation*}\label{eq:covariancelonger}
E_{Q}[\mathcal{Y}_{t}(H)\mathcal{Y}_{s}(G)]=\chi(\rho)\int_{\mathbb{R}}H(u)G(u)du
\end{equation*}
for every $s,t\geq{0}$ and $H$, $G$ in $\mathcal{H}_{k}$. Then, $(Q^{\gamma}_{N})_{N\geq{1}}$ converges weakly to $Q$. \label{th:flu2}
\end{theorem}

The main problem to overcome when showing last result is the Boltzmann-Gibbs Principle, which we can prove for $\gamma<1/3$ using a multi-scale
argument as for the ASEP in \cite{G.}.

\begin{theorem}{(Boltzmann-Gibbs Principle)} \label{th:bg}

Fix $\gamma<1/3$. For every $t>0$ and $H\in{\mathcal{S}(\mathbb{R})}$,
\begin{equation*}
\lim_{N\rightarrow{\infty}}\mathbb{E}_{\nu_{\rho}}^{\gamma}\Big[\int_{0}^{t}
\frac{N^{\gamma}}{\sqrt{N}}\sum_{x\in{\mathbb{Z}}}H\Big(\frac{x}{N}\Big)V_{g}(\eta_{s}(x))ds\Big]^{2}=0,
\end{equation*}
where
\begin{equation*}
V_{g}(\eta(x))=g(\eta(x))-\phi(\rho)-\phi'(\rho)[\eta(x)-\rho],
\end{equation*}
and $\phi(\rho)=E_{\nu_{\rho}}[g(\eta(0))]$.
\end{theorem}

Now we define the current of particles across a characteristic. Let $J^{N,\gamma}_{v^{x}_{t}}(tN)$ be the current through the bond
$[v^{x}_{t},v^{x}_{t}+1]$ (where $v^{x}_{t}=x+[\phi'(\rho)tN^{1+\gamma}]$) defined as the number of particles that jump from $v^{x}_{t}$ to
$v^{x}_{t}+1$, from time $0$ to $tN^{1+\gamma}$:
\begin{equation*}
J^{N,\gamma}_{v^{x}_{t}}(tN)=\sum_{y\geq{1}}\Big(\eta_{t}(y+v^{x}_{t})-\eta_{0}(y+x)\Big).
\end{equation*}
As a consequence of last result, it holds that:
\begin{proposition} \label{current dependence higher}
Fix $t\geq{0}$, a site $x\in{\mathbb{Z}}$ and $\gamma<1/3$. Then,
\begin{equation*}
\lim_{N\rightarrow{+\infty}}\mathbb{E}_{\nu_{\rho}}^{\gamma}\Big[\frac{\bar{J}^{N,\gamma}_{v^{x}_{t}}(tN)}{\sqrt{N}}\Big]^{2}=0.
\end{equation*}
\end{proposition}

\section{Density Fluctuations for the Hyperbolic Scaling}

\subsection{Equilibrium Fluctuations}
Fix a positive integer $k$, denote by $\mathfrak{A}$ the operator $\phi'(\rho)\nabla$ defined on a domain of $L^{2}(\mathbb{R})$ and by
$\{T_{t}, t\geq0\}$ its semigroup. The theorem follows as long as we show that $(Q_{N})_{N\geq{1}}$ is tight and characterize the limiting
measure $Q$.

Fix $H\in{\mathcal{S}(\mathbb{R})}$, then
\begin{equation*} \label{eq:martingale M}
M^{N,H}_{t}=\mathcal{Y}^{N}_{t}(H)-\mathcal{Y}^{N}_{0}(H)-\int^{t}_{0}\frac{1}{\sqrt{N}}\sum_{x\in\mathbb{Z}}\nabla^{N}{H\Big(\frac{x}{N}\Big)}g(\eta_{s}(x))ds
\end{equation*}
is a martingale with respect to the natural filtration with quadratic variation given by
\begin{equation*} \label{eq:quadraticvariation}
\int^{t}_{0}\frac{1}{N^2}\sum_{x\in\mathbb{Z}}\Big(\nabla^{N}{H\Big(\frac{x}{N}\Big)}\Big)^2\Big[g(\eta_{s}(x))+g(\eta_{s}(x+1))\Big]ds,
\end{equation*}
where $\nabla^{N} H$ denotes the discrete derivative of $H$. The integral part of the martingale can be written as
\begin{equation*}
\int^{t}_{0}\frac{1}{\sqrt{N}}\sum_{x\in\mathbb{Z}}\nabla^{N}{H\Big(\frac{x}{N}\Big)} \Big[g(\eta_{s}(x))-\phi(\rho)\Big]ds.
\end{equation*}
by using the fact that $\sum_{x\in{\mathbb{Z}}}\nabla^{N}{H(\frac{x}{N})}=0$. The following result allows to replace $g(\eta_{s}(x))-\phi(\rho)$
by $\phi'(\rho)[\eta_{s}(x)-\rho]$ and allows to recover the density fluctuation field inside the integral part of the martingale.

\begin{theorem}{(Boltzmann-Gibbs Principle)}

For every $H\in{\mathcal{S}(\mathbb{R})}$ and every $t>0$,
\begin{equation*}
\lim_{N\rightarrow{\infty}}\mathbb{E}_{\nu_{\rho}}\Big[\Big(\int_{0}^{t}
\frac{1}{\sqrt{N}}\sum_{x\in{\mathbb{Z}}}H\Big(\frac{x}{N}\Big)V_{g}(\eta_{s}(x))ds\Big)^2\Big]=0.
\end{equation*}
\end{theorem}

The proof of last result follows the same lines as for the Symmetric Zero-Range Process in \cite{K.L.} and for that reason we have omitted it.
For the same reason the following results are just stated but their proofs follow the same lines as for the ASEP in \cite{G.}:
$(Q_{N})_{N\geq{1}}$ is a tight sequence, the limiting measure $Q$ is supported on fields $\mathcal{Y}_{\cdot}$ such that for a fixed time $t$
and a test function $H$, $\mathcal{Y}_{t}(H)=\mathcal{Y}_{0}(T_{t}H)$ where $T_{t}H(u)=H(u+\phi'(\rho)t)$ and $\mathcal{Y}_{0}$ is a Gaussian
field with covariance given by $E_{Q}(\mathcal{Y}_{0}(G)\mathcal{Y}_{0}(H))=\chi(\rho)<G,H>$.

\subsection{Central Limit Theorem for the Current over a fixed bond}

Now we give a sketch of the proof of Theorem \ref{th:cltcurrent} in which we need to show the convergence of finite dimensional distributions of
$Z_{t}^{N}/\sqrt{\chi(\rho)\phi'(\rho)}$ to those of Brownian motion together with tightness.

We start by the convergence of finite dimensional distributions, namely, we show that for every $k\geq{1}$ and every
$0\leq{t_{1}}<{t_{2}}<..<t_{k}$, $(Z_{t_{1}}^{N},..,Z_{t_{k}}^{N})$ converges in law to a Gaussian vector $(Z_{t_{1}},..,Z_{t_{k}})$ with mean
zero and covariance given by $E_{Q}[Z_{t}Z_{s}]=\chi(\rho)\phi'(\rho)s$ provided $s\leq{t}$.

 Recall that
$J^N_{-1,0}(tN)$ is defined as the total number of jumps from the site $-1$ to $0$ during the time interval $[0,tN]$. Since
\begin{equation*}
J^N_{-1,0}(tN)=\sum_{x\geq{0}}\Big(\eta_{t}(x)-\eta_{0}(x)\Big),
\end{equation*}
the current can be written in terms of the density fluctuation field evaluated on $H_{0}$, the Heaviside function $H_{0}(u)=1_{[0,\infty)}(u)$:
\begin{equation*}
\frac{1}{\sqrt{N}}\Big\{J^N_{-1,0}(tN)-\mathbb{E}_{\nu_{\rho}}[J^N_{-1,0}(tN)]\Big\}=\mathcal{Y}_{t}^{N}(H_{0})-\mathcal{Y}_{0}^{N}(H_{0}),
\end{equation*}
Approximating $H_{0}$ by $(G_{n})_{n\geq{1}}$ such that $G_{n}(u)=(1-\frac{u}{n})^{+}1_{[0,\infty)}(u)$, then
\begin{proposition} \label{prop:1}
For every $t\geq{0}$,
\begin{equation*}
\lim_{n\rightarrow{+\infty}}\mathbb{E}_{\nu_{\rho}}\Big[\frac{\bar{J}^N_{-1,0}(tN)}{\sqrt{N}}-(\mathcal{Y}_{t}^{N}(G_{n})-\mathcal{Y}_{0}^{N}(G_{n}))\Big]^{2}=0
\end{equation*}
uniformly in $N$.
\end{proposition}
The convergence of finite dimensional distributions is an easy consequence of last result together with Theorem \ref{th:flu1}, see \cite{J.L.}.

Now, it remains to prove that the distributions of $Z_{t}^{N}/\sqrt{\chi(\rho)\phi'(\rho)}$ are tight. For that, we can use the same argument as
in Theorem 2.3 of \cite{G.} that relies on the use of Theorem 2.1 of \cite{S.} with the definition of weakly positive associated increments
given in \cite{S.2}. One can follow the same arguments as those of Theorem 2 of \cite{K.} to show that $J_{-1,0}(t)$ has weakly positiive
associated increments with the definition in \cite{S.2}, see \cite{G.}. In order to conclude the proof it remains to note that
\begin{equation*}
\lim_{N\rightarrow{+\infty}}\frac{1}{t}E_{\nu_{\rho}}[J_{-1,0}(t)]=\sigma^2,
\end{equation*}
which follows by Theorem 3 of \cite{K.}.

\section{Density Fluctuations for a longer time scale}
Fix a positive integer $k$ and let $U_{t}^{N}H(u)=H(u-\phi'(\rho)tN^{\gamma})$. Recall the definition of $(Q^\gamma_{N})_{N\geq{1}}$ and note
that following the same computations as in \cite{G.} it is easy to show that the sequence is tight. Now we compute the limit field, by fixing
$H\in{\mathcal{S}(\mathbb{R})}$ such that
\begin{equation*}
M_{t}^{N,H}=\mathcal{Y}_{t}^{N,\gamma}(H)-\mathcal{Y}_{0}^{N,\gamma}(H)-\int_{0}^{t}\frac{N^{\gamma}}{\sqrt{N}}\sum_{x\in{\mathbb{Z}}}\nabla^{N}
U_{s}^{N}H\Big(\frac{x}{N}\Big)V_{g}(\eta_{s}(x))
\end{equation*}
is a martingale and whose quadratic variation is given by
\begin{equation*} \label{quadraticlonger}
\int^{t}_{0}\frac{N^{\gamma}}{N^2}\sum_{x\in\mathbb{Z}}
\Big(\nabla^{N}{U_{s}^{N}H\Big(\frac{x}{N}\Big)}\Big)^2\Big[g(\eta_{s}(x))+g(\eta_{s}(x+1))\Big]ds.
\end{equation*}

If $\gamma<1$, $M_{t}^{N,H}$ vanishes in $L^{2}(\mathbb{P}_{\nu_\rho}^{\gamma})$ as $N\rightarrow{+\infty}$. Using the Botzmann-Gibbs Principle,
whose proof is sketched in the next section, the integral part of the martingale $M_{t}^{N,H}$ vanishes in
$L^{2}(\mathbb{P}_{\nu_\rho}^{\gamma})$ as $N\rightarrow{+\infty}$ which in turn implies that if $Q$ is one limiting point of $(Q_{N})_{N}$, the
limit density fluctuation field satisfies $\mathcal{Y}_{t}(H)=\mathcal{Y}_{0}(H)$, where $\mathcal{Y}_{0}$ is a Gaussian field with covariance
given by $E_{Q}(\mathcal{Y}_{0}(G)\mathcal{Y}_{0}(H))=\chi(\rho)<G,H>$.

\section{Boltzmann-Gibbs Principle}

In this section we prove Theorem \ref{th:bg}. Since we are going to follow the same steps as in Theorem (2.6) of \cite{G.} we just remark the
fundamental differences between the proofs.

To start fix an integer $K$ and a test function $H\in \mathcal{S}(\mathbb{\mathbb{R}})$. We divide $\mathbb{Z}$ in non overlapping intervals of
length $K$, denoted by $\{I_{j},j\geq{1}\}$ and by summing and subtracting $H\Big(\frac{y_{j}}{N}\Big)$, where $y_{j}$ is some point of $I_{j}$,
we can bound the expectation appearing in the statement of the Theorem by
\begin{equation*}
2\mathbb{E}_{\nu_{\rho}}^{\gamma}\Big[\int_{0}^{t}\frac{N^{\gamma}}{\sqrt{N}}\sum_{j\geq{1}}
\sum_{x\in{I}_{j}}\Big[H\Big(\frac{x}{N}\Big)-H\Big(\frac{y_{j}}{N}\Big)\Big]V_{g}(\eta_{s}(x))ds\Big]^{2}
\end{equation*}
\begin{equation*}\label{eq:bg3}
+2\mathbb{E}_{\nu_{\rho}}^{\gamma}\Big[\int_{0}^{t}\frac{N^{\gamma}}{\sqrt{N}}\sum_{j\geq{1}}
H\Big(\frac{y_{j}}{N}\Big)\sum_{x\in{I_{j}}}V_{g}(\eta_{s}(x))ds\Big]^{2}.
\end{equation*}

The first expectation is easily handled, since by Schwarz inequality and the invariance of $\nu_{\rho}$ it can be bounded by
$Ct^2N^{2\gamma}||H'||_{2}^{2}\Big(\frac{K}{N}\Big)^2$ and vanishes as long as $K N^{\gamma-1}\rightarrow{0}$  when $N\rightarrow{+\infty}$.

In order to treat the remaining expectation we bound it from above by
\begin{equation} \label{eq:bg4}
2\mathbb{E}_{\nu_{\rho}}^{\gamma}\Big[\int_{0}^{t}\frac{N^{\gamma}}{\sqrt{N}}\sum_{j\geq{1}}H\Big(\frac{y_{j}}{N}\Big)
V_{1,j,g}(\eta_{s})ds\Big]^{2}
+2\mathbb{E}_{\nu_{\rho}}^{\gamma}\Big[\int_{0}^{t}\frac{N^{\gamma}}{\sqrt{N}}\sum_{j\geq{1}}H\Big(\frac{y_{j}}{N}\Big)E\Big(\sum_{x\in{I_{j}}}
V_{g}(\eta_{s}(x))\Big|M_{j}\Big)ds\Big]^{2}
\end{equation}
where
\begin{equation*}
V_{1,j,g}(\eta)=\sum_{x\in{I_{j}}} V_{g}(\eta(x))-E\Big(\sum_{x\in{I_{j}}} V_{g}(\eta(x))\Big|M_{j}\Big).
\end{equation*}
and $M_{j}=\sigma\Big(\sum_{x\in{I_{j}}}\eta(x)\Big)$.
\begin{lemma} \label{lm:lengthK}
For every $H\in \mathcal{S}(\mathbb{R})$ and every $t>0$, if $K^{2}N^{\gamma-1}\rightarrow{0}$ as $N\rightarrow{+\infty}$, then
\begin{equation*}
\lim_{N\rightarrow{\infty}}\mathbb{E}_{\nu_{\rho}}^{\gamma}\Big[\int_{0}^{t}\frac{N^{\gamma}}{\sqrt{N}}\sum_{j\geq{1}}H\Big(\frac{y_{j}}{N}\Big)
V_{1,j,g}(\eta_{s})ds\Big]^{2}=0.
\end{equation*}
\end{lemma}
\begin{proof}
By Proposition A1.6.1 of \cite{K.L.} and by the variational formula for the $H_{-1}$-norm the expectation above is bounded by
\begin{equation*}
Ct \sum_{j\geq{1}}\sup_{h \in L^{2}(\nu_{\rho})}\Big\{2\int
\frac{N^{\gamma}}{\sqrt{N}}H\Big(\frac{y_{j}}{N}\Big)V_{1,j,g}(\eta)h(\eta)\nu_{\rho}(d\eta)-N^{1+\gamma}<h,-\mathcal{L}^{S}_{I_{j}}h>_{\rho}\Big\},
\end{equation*}
where $\mathcal{L}^{S}$ is the Symmetric dynamics restricted to the set $I_{j}$, namely:
\begin{equation*}
\mathcal{L}^{S}_{I_{j}}f(\eta)=\sum_{\substack{x,y\in{I_{j}}\\|x-y|=1}}\frac{1}{2}1_{\{\eta(x)\geq{1}\}}[f(\eta^{x,y})-f(\eta)].
\end{equation*}
where
\[\eta^{x,y}(z)=
\begin{cases}
\eta(z), & \mbox{if $z\neq{x,y}$}\\
\eta(x)-1, & \mbox{if $z=x$}\\
\eta(y)+1, & \mbox{if $z=y$}
\end{cases}.
\]
\\

For each $j$ and $A_{j}$ a positive constant, it holds that
\begin{equation*}
\int V_{1,j,g}(\eta)h(\eta)\nu_{\rho}(d\eta)\leq{\frac{1}{2
A_{j}}<V_{1,j,g},(-\mathcal{L}^{S}_{I_{j}})^{-1}V_{1,j,g}>_{\rho}+\frac{A_{j}}{2}<h,-\mathcal{L}^{S}_{I_{j}}h>_{\rho}},
\end{equation*}
and taking $A_{j}=N^{3/2}\Big(|H(\frac{y_{j}}{N})|\Big)^{-1}$, the whole expectation becomes bounded by
\begin{equation*}
Ct\sum_{j\geq{1}}\frac{N^{\gamma}}{N^{2}}{H^{2}\Big(\frac{y_{j}}{N}\Big)}<V_{1,j,g},(-\mathcal{L}^{S}_{I_{j}})^{-1}V_{1,j,g}>_{\rho}.
\end{equation*}
 By the spectral gap inequality for the Symmetric Zero-Range process (see
\cite{L.S.V.}) last expression can be bounded by
\begin{equation*}
Ct\sum_{j\geq{1}}\frac{N^{\gamma}}{N^{2}}{H^{2}\Big(\frac{y_{j}}{N}\Big)}(K+1)^{2}Var(V_{1,j,g},\nu_{\rho}).
\end{equation*}
The proof of the Lemma ends if we show that $Var(V_{1,j,g},\nu_{\rho})\leq{K C}$, since it implies that the expectation in the statement of the
lemma to be bounded by $C t\frac{N^{\gamma}}{N}(K+1)^{2}||H||_{2}^{2}$ and vanishes as long as $K^{2}N^{\gamma-1}\rightarrow{0}$ when
$N\rightarrow{+\infty}$.
\end{proof}

\begin{remark}
Here we show that $Var(V_{1,j,g},\nu_{\rho})\leq{K C}$. Since $Var(V_{1,j,g},\nu_{\rho})\leq{E_{\nu_{\rho}}[V_{1,j,g}]^2}$ and by the definition
of $V_{1,j,g}$  we have that
\begin{equation*}
Var(V_{1,j,g},\nu_{\rho})\leq{E_{\nu_{\rho}}\Big[\sum_{x\in{I_{j}}}V_{g}(\eta(x))-
E_{\nu_{\rho}}[\sum_{x\in{I_{j}}}V_{g}(\eta(x))|M_{j}]\Big]^2}.
\end{equation*}
By the definition of $V_{g}(\eta)$ last expression can be written as
\begin{equation*}
E_{\nu_{\rho}}\Big[\sum_{x\in{I_{j}}}\Big(g(\eta(x))-\phi(\rho)-\phi'(\rho)[\eta(x)-\rho]\Big)-\sum_{x\in{I_{j}}}\phi_{j}(\rho)
-K\phi(\rho)-\sum_{x\in{I_{j}}}\phi'(\rho)[\eta_{j}^{K}-\rho]\Big]^2,
\end{equation*}
where $\phi_{j}(\rho)=E_{\nu_{\rho}}[g(\eta)|M_{j}]$ and $\eta_{j}^{K}=\frac{1}{K}\sum_{x\in{I_{j}}}\eta(x)$. On the other hand, by summing and
subtracting $\phi(\eta_{j}^{K})=E_{\nu_{\eta_{j}^{K}}}[g(\eta)]$, where $\nu_{\eta^{K}_{j}}$ is the Bernoulli measure with density
$\eta^{K}_{j}$, last expression can be bounded by
\begin{equation*}
2E_{\nu_{\rho}}\Big[\sum_{x\in{I_{j}}}\Big(g(\eta(x))-\phi(\rho)\Big)\Big]^2
+2E_{\nu_{\rho}}\Big[\sum_{x\in{I_{j}}}\phi'(\rho)[\eta(x)-\rho]\Big]^2
\end{equation*}
\begin{equation*}
2E_{\nu_{\rho}}\Big[\sum_{x\in{I_{j}}}\Big (\phi_{j}(\rho)-\phi(\eta_{j}^{K})\Big)\Big]^2
+2E_{\nu_{\rho}}\Big[\sum_{x\in{I_{j}}}\Big(\phi(\eta_{j}^{K})-\phi(\rho)-\phi'(\rho)[\eta_{j}^{K}-\rho]\Big)\Big]^2.
\end{equation*}
Now we treat each expectation separately.

For the first and the second one, since $(\eta(x))_{x}$ are independent under $\nu_{\rho}$, it is easy to show that
\begin{equation*}
E_{\nu_{\rho}}\Big[\sum_{x\in{I_{j}}}\Big(g(\eta(x))-\phi(\rho)\Big)\Big]^2\leq{C K Var(g,\nu_{\rho})}
\end{equation*}
and
\begin{equation*}
E_{\nu_{\rho}}\Big[\sum_{x\in{I_{j}}}\phi'(\rho)[\eta(x)-\rho]\Big]^2\leq{C K Var(\eta(0),\nu_{\rho})}.
\end{equation*}
On the other hand, to treat the third expectation one can use the equivalence of ensembles (see Corollary A2.1.7 of \cite{K.L.}) which
guarantees that $|\phi_{j}(\rho)-\phi(\eta_{j}^{K})|\leq{\frac{C(g)}{K}}$ while for the last one, one can use Taylor expansion to have
\begin{equation*}
E_{\nu_{\rho}}\Big[\phi(\eta_{j}^{K})-\phi(\rho)-\phi'(\rho)(\eta_{j}^{K}-\rho)\Big]^2\sim
{E_{\nu_{\rho}}\Big[\eta_{j}^{K}-\rho\Big]^4=O(K^{-2})}.
\end{equation*}
Putting these arguments all together one gets to the bound $KC$.
\end{remark}

To conclude the proof it remains to bound the expectation on the right hand side of (\ref{eq:bg4}). For that, fix an integer $L$ and take
disjoint intervals of length $M=LK$, denoted by $\{\tilde{I}_{l},l\geq{1}\}$ and write it as:
\begin{equation*}
\mathbb{E}_{\nu_{\rho}}^{\gamma}\Big[\int_{0}^{t}\frac{N^{\gamma}}{\sqrt{N}}\sum_{l\geq{1}}\sum_{j\in{\tilde{I}_{l}}}H\Big(\frac{y_{j}}{N}\Big)
E\Big(\sum_{x\in{I_{j}}}V_{g}(\eta_{s}(x))\Big|M_{j}\Big)ds\Big]^{2}.
\end{equation*}
By summing and subtracting $H\Big(\frac{z_{l}}{N}\Big)$, where $z_{l}$ denotes one point of the interval $\tilde{I}_{l}$, last expectation can
be bounded by
\begin{equation*}
2\mathbb{E}_{\nu_{\rho}}^{\gamma}\Big[\int_{0}^{t}\frac{N^{\gamma}}{\sqrt{N}}\sum_{l\geq{1}}\sum_{j\in{\tilde{I}_{l}}}
\Big[H\Big(\frac{y_{j}}{N}\Big)-H\Big(\frac{z_{l}}{N}\Big)\Big] E\Big(\sum_{x\in{I_{j}}}V_{g}(\eta_{s}(x))\Big|M_{j}\Big)ds\Big]^{2}
\end{equation*}
\begin{equation*} \label{eq:bg7}
+2\mathbb{E}_{\nu_{\rho}}^{\gamma}\Big[\int_{0}^{t}\frac{N^{\gamma}}{\sqrt{N}}\sum_{l\geq{1}}H\Big(\frac{z_{l}}{N}\Big)
\sum_{j\in{\tilde{I}_{l}}}E\Big(\sum_{x\in{I_{j}}}V_{g}(\eta_{s}(x))\Big|M_{j}\Big)ds\Big]^{2}.
\end{equation*}
Following the same arguments as above it is easy to show that the first expectation vanishes if $L^{2}KN^{2\gamma-2}\rightarrow{0}$ as
$N\rightarrow{+\infty}$. For the second one, sum and subtract $E\Big(\sum_{x\in{\tilde{I}_{l}}}V_{g}(\eta(x))\Big|\tilde{M}_{l}\Big)$ where
$\tilde{M}_{l}=\sigma\Big(\sum_{x\in{\tilde{I}_{l}}}\eta(x)\Big)$ and bound it by
\begin{equation*}\label{eq:bg8}
2\mathbb{E}_{\nu_{\rho}}^{\gamma}\Big[\int_{0}^{t}\frac{N^{\gamma}}{\sqrt{N}}\sum_{l\geq{1}}
H\Big(\frac{z_{l}}{N}\Big)V_{2,l,g}(\eta_{s})ds\Big]^{2}
+2\mathbb{E}_{\nu_{\rho}}^{\gamma}\Big[\int_{0}^{t}\frac{N^{\gamma}}{\sqrt{N}}\sum_{l\geq{1}}H\Big(\frac{z_{l}}{N}\Big)
E\Big(\sum_{x\in{\tilde{I}_{l}}}V_{g}(\eta_{s}(x))\Big|\tilde{M}_{l}\Big)ds\Big]^{2},
\end{equation*}
 where
\begin{equation*}
V_{2,l,g}(\eta)=\sum_{j\in{\tilde{I}_{l}}}E\Big(\sum_{x\in{I_{j}}}
V_{g}(\eta(x))\Big|M_{j}\Big)-E\Big(\sum_{x\in{\tilde{I}_{l}}}V_{g}(\eta(x))\Big|\tilde{M}_{l}\Big).
\end{equation*}

\begin{lemma} \label{lm:lengthM}
For every $H\in \mathcal{S}(\mathbb{R})$ and every $t>0$, if $L^{2}KN^{\gamma-1}\rightarrow{0}$ as $N\rightarrow{+\infty}$, then
\begin{equation*}
\lim_{N\rightarrow{\infty}}\mathbb{E}_{\nu_{\rho}}^{\gamma}\Big[\int_{0}^{t}\frac{N^{\gamma}}{\sqrt{N}}\sum_{l\geq{1}}H\Big(\frac{z_{l}}{N}\Big)
V_{2,l,g}(\eta_{s})ds\Big]^{2}=0.
\end{equation*}
\end{lemma}
\begin{proof}
Following the proof of Lemma (\ref{lm:lengthK}), the expectation becomes bounded by
\begin{equation*}
Ct \sum_{l\geq{1}}\sup_{h \in L^{2}(\nu_{\rho})}\Big\{2 \int \frac{N^{\gamma}}{\sqrt{N}}
H\Big(\frac{z_{l}}{N}\Big)V_{2,l,g}(\eta)h(\eta)\nu_{\rho}(d\eta)-N^{1+\gamma}<h,-\mathcal{L}^{S}_{\tilde{I}_{l}}h>_{\rho}\Big\}.
\end{equation*}
 Using an appropriate $A_{l}$ and the spectral gap inequality, last expression is bounded by
\begin{equation*}
Ct\sum_{l\geq{1}}\frac{N^{\gamma}}{N^{2}}H^{2}\Big(\frac{z_{l}}{N}\Big)(M+1)^{2}Var(V_{2,l,g},\nu_{\rho}).
\end{equation*}
Now, the proof ends as long as $Var(V_{2,l,g},\nu_{\rho})\leq{L C}$, which is proved below.
\end{proof}

\begin{remark}
Here we show that $Var(V_{2,l,g},\nu_{\rho})\leq{L C}$. Since $Var(V_{2,l,g},\nu_{\rho})\leq{E_{\nu_{\rho}}[V_{2,l,g}]^2}$ and by the definition
of $V_{2,l,g}$ we have that
\begin{equation*}
Var(V_{2,l,g},\nu_{\rho})\leq{E_{\nu_{\rho}}\Big[\sum_{j\in{\tilde{I}_{l}}}E\Big(\sum_{x\in{I_{j}}}
V_{g}(\eta(x))\Big|M_{j}\Big)-E\Big(\sum_{x\in{\tilde{I}_{l}}}V_{g}(\eta(x))\Big|\tilde{M}_{l}\Big)\Big]^2}.
\end{equation*}
By the definition of $V_{g}(\eta)$ and the notation introduced above, one can write last expression as
\begin{equation*}
E_{\nu_{\rho}}\Big[\sum_{j\in{\tilde{I}_{l}}}\Big(K\phi_{j}(\rho)-K\phi(\rho)-\phi'(\rho)K[\eta_{j}^{K}-\rho]\Big)-
M\phi_{l}(\rho)-M\phi(\rho)-M\phi'(\rho)[\eta_{l}^{M}-\rho]\Big]^2,
\end{equation*}
where $\phi_{l}(\rho)=E_{\nu_{\rho}}[g(\eta)|M_{l}]$ and $\eta_{l}^{M}=\frac{1}{M}\sum_{x\in{I_{l}}}(\eta(x)-\rho)$. Last expression can be
written as
\begin{equation*}
E_{\nu_{\rho}}\Big[M\Big\{\frac{1}{M}\sum_{j\in{\tilde{I}_{l}}}\Big(K\phi_{j}(\rho)-K\phi(\rho)-\phi'(\rho)K[\eta_{j}^{K}-\rho]\Big)-
\phi_{l}(\rho)-\phi(\rho)-\phi'(\rho)[\eta_{l}^{M}-\rho]\Big\}\Big]^2
\end{equation*}
\begin{equation*}
=E_{\nu_{\rho}}\Big[M\Big\{\frac{1}{L}\sum_{j\in{\tilde{I}_{l}}}\Big(\phi_{j}(\rho)-\phi'(\rho)[\eta_{j}^{K}-\rho]-
\phi_{l}(\rho)-\phi'(\rho)[\eta_{l}^{M}-\rho]\Big)\Big\}\Big]^2
\end{equation*}
\begin{equation*}
=\frac{M^{2}}{L}E_{\nu_{\rho}}\Big[\frac{1}{\sqrt{L}}\sum_{j\in{\tilde{I}_{l}}}\Big(\phi_{j}(\rho)-\phi'(\rho)[\eta_{j}^{K}-\rho]-
\phi_{l}(\rho)-\phi'(\rho)[\eta_{l}^{M}-\rho]\Big)\Big]^2.
\end{equation*}

By the independence of the random variables $(\eta(x))_{x}$ under $\nu_{\rho}$ and the Central Limit Theorem, last expectation is of order
\begin{equation*}
E_{\nu_{\rho}}\Big[\phi_{j}(\rho)-\phi(\rho)-\phi'(\rho)[\eta_{j}^{K}-\rho]\Big]^2,
\end{equation*}
 which we can bound by
\begin{equation*}
2E_{\nu_{\rho}}\Big[\phi_{j}(\rho)-\phi(\eta_{j}^{K})\Big]^2+2
E_{\nu_{\rho}}\Big[\phi(\eta_{j}^{K})-\phi(\rho)-\phi'(\rho)[\eta_{j}^{K}-\rho]\Big]^2.
\end{equation*}
By the equivalence of ensembles the expectation on the left hand side is bounded by $K^{-2}$. For the other, use Taylor expansion to have
\begin{equation*}
E_{\nu_{\rho}}\Big[\phi(\eta_{j}^{K})-\phi(\rho)-\phi'(\rho)[\eta_{j}^{K}-\rho]\Big]^2\sim{ E_{\nu_{\rho}}[\eta_{j}^{K}-\rho]^4=O(K^{-2})}
\end{equation*}
This finishes the proof of the remark.

\end{remark}

\vspace{0,2cm}

\textbf{The proof of Boltzmann-Gibbs Principle}

\vspace{0,2cm}

Following the same arguments as before, take $n$ sufficiently big for which in the n-th step of the proof we have intervals, denoted by
$\{I^{n}_{p},p\geq{1\geq}\}$ of length $K_{n}=N^{1-\gamma}$. At this stage it remains to bound:
\begin{equation*}
\mathbb{E}_{\nu_{\rho}}^{\gamma}\Big[\int_{0}^{t}\frac{N^{\gamma}}{\sqrt{N}}\sum_{p\geq{1}}H\Big(\frac{z_{p}}{N}\Big)
E_{\nu_{\rho}}\Big(\sum_{x\in{{I}^{n}_{p}}}V_{g}(\eta_{s}(x))\Big|{M}_{p}^{n}\Big)ds\Big]^{2},
\end{equation*}
where for each $p$, $z_{p}$ is one point of the interval $I^{n}_{p}$ and ${M}_{p}^{n}=\sigma\Big(\sum_{x\in{{I}^{n}_{p}}}\eta(x)\Big)$.

Since $\nu_{\rho}$ is an invariant product measure, last expectation can be bounded by
\begin{equation} \label{eq:necessaria BG fim}
t^{2}\frac{N^{2\gamma}}{N}\sum_{p\geq{1}}\Big(H\Big(\frac{z_{p}}{N}\Big)\Big)^{2}{E}_{\nu_{\rho}}\Big(
E_{\nu_{\rho}}\Big(\sum_{x\in{{I}^{n}_{p}}}V_{g}(\eta(x))\Big|{M}_{p}^{n}\Big)\Big)^{2}.
\end{equation}

\begin{remark}
Here we show that $E_{\nu_{\rho}}\Big(E_{\nu_{\rho}}\Big(\sum_{x\in{{I}^{n}_{p}}}V_{g}(\eta(x))\Big|M_{p}^{n}\Big)\Big)^{2}=O(1)$.

By the definition of $V_{g}$, the expectation above is equal to
\begin{equation*}
E_{\nu_{\rho}}\Big(E_{\nu_{\rho}}\Big(\sum_{x\in{{I}^{n}_{p}}}\Big(\phi_{K_{n}}(\rho)-\phi(\rho)-\phi'(\rho)[\eta_{n}^{K_{n}}-\rho]\Big)\Big)^{2}
\end{equation*}
and bounded from above by
\begin{equation*}
2E_{\nu_{\rho}}\Big(E_{\nu_{\rho}}\Big(\sum_{x\in{{I}^{n}_{p}}}\Big(\phi_{K_{n}}(\rho)-\phi(\eta_{n}^{K_{n}})\Big)\Big)^{2}+
2E_{\nu_{\rho}}\Big(E_{\nu_{\rho}}\Big(\sum_{x\in{{I}^{n}_{p}}}\Big(\phi(\eta_{n}^{K_{n}})-\phi(\rho)-\phi'(\rho)[\eta_{n}^
{K_{n}}-\rho]\Big)\Big)\Big)^{2},
\end{equation*}
where $\phi_{K_{n}}(\rho)=E_{\nu_{\rho}}[g(\eta)|M_{p}^{n}]$ and $\eta_{n}^{K_{n}}=\frac{1}{K_{n}}\sum_{x\in{I_{p}^{n}}}\eta(x)$. Now the result
follows if one applies equivalence of ensembles to the expectation on the left hand side and Taylor expansion to the expectation on the right
hand side.
\end{remark}
This implies (\ref{eq:necessaria BG fim}) to be bounded by $\frac{N^{2\gamma}}{K_{n}}$, which vanishes as $N\rightarrow{+\infty}$ since
$\gamma<1/3$.
\begin{remark}
 Here we give an application of the Boltzmann-Gibbs Principle for a linear functional associated to the
one-dimensional Symmetric Zero-Range process, in the diffusive scaling. Consider a Markov process $\eta_{tN^2}$ with generator given by
\begin{equation*}
\mathcal{L}^{S}f(\eta)=\sum_{\substack{x,y\in{\mathbb{Z}}\\|x-y|=1}}\frac{1}{2}1_{\{\eta(x)\geq{1}\}}[f(\eta^{x,y})-f(\eta)],
\end{equation*}
with $\eta^{x,y}$ as defined in the proof of Lemma \ref{lm:lengthK}.
 If one repeats the same steps as done in the proof of Theorem \ref{th:bg} it is easy to show that:
\begin{corollary}
Fix $t>0$ and $\beta<1/2$, then
\begin{equation*}
\lim_{N\rightarrow{\infty}}\mathbb{E}_{\nu_{\rho}}\Big[{N^{\beta}}\int_{0}^{t}
\frac{1}{\sqrt{N}}\sum_{x\in{\mathbb{Z}}}H\Big(\frac{x}{N}\Big)V_{g}(\eta_{s}(x))ds\Big]^{2}=0.
\end{equation*}
\end{corollary}
So, in order to observe fluctuations for this field one has to take $\beta\geq{1/2}$.
\end{remark}

\subsection{Current through the characteristics speed}

As in the hyperbolic scaling, Proposition (\ref{current dependence higher}) is a consequence of:
\begin{proposition}
For every $t\geq{0}$ and $\gamma<1/3$:
\begin{equation*}
\lim_{n\rightarrow{+\infty}}\mathbb{E}_{\nu_{\rho}}^{\gamma}\Big[\frac{\bar{J}^{N,\gamma}_{v^{x}_{t}}(tN)}
{\sqrt{N}}-(\mathcal{Y}_{t}^{N,\gamma}(G_{n})-\mathcal{Y}_{0}^{N,\gamma}(G_{n}))\Big]^{2}=0,
\end{equation*}
uniformly over $N$.
\end{proposition}
The proof of this result follows the same lines as the proof of Proposition 9.4 in \cite{G.} and for that reason we have omitted it.

\end{document}